\documentclass[12pt]{amsart}
\usepackage{amscd, amssymb}
\begin{document}

\def\supp{\operatorname{supp}}
\def\Ind{\operatorname{Ind}}
\def\id{\textnormal{id}}
\def\Aut{\operatorname{Aut}}
\def\range{\operatorname{range}}
\def\sp{\operatorname{span}}
\def\clsp{\overline{\operatorname{span}}}
\def\Prim{\operatorname{Prim}}
\def\dashind{\operatorname{\!-Ind}}
\def\rt{\textnormal{rt}}
\def\lt{\textnormal{lt}}
\def\tr{\operatorname{tr}}

\def\H{\mathcal{H}}
\def\L{\mathcal{L}}
\def\K{\mathcal{K}}
\def\B{B}
\def\C{\mathbb{C}}

\newtheorem{thm}{Theorem}
\newtheorem{cor}[thm]{Corollary}
\newtheorem{lemma}[thm]{Lemma}
\newtheorem{prop}[thm]{Proposition}
\newtheorem{thm1}{Theorem}

\theoremstyle{definition}
\newtheorem{definition}[thm]{Definition}

\theoremstyle{remark}
\newtheorem{remark}[thm]{Remark}
\newtheorem{example}[thm]{Example}
\newtheorem{remarks}[thm]{Remarks}
\newtheorem{claim}[thm]{Claim}

\title[Regularity of induced representations]{\boldmath  Regularity
of induced representations and\\ a theorem of Quigg and Spielberg }

\author[Astrid an Huef]{Astrid an Huef}
\address{Department of Mathematics and Computer Science\\
University of Denver\\
2360 S. Gaylord St.\\
Denver, CO 80208-0189\\
USA}
\email{astrid@cs.du.edu}

\author[Iain Raeburn]{Iain Raeburn}
\address{Department of Mathematics\\ University of Newcastle\\
NSW 2308\\ Australia}
\email{iain@frey.newcastle.edu.au}

\thanks{This research was supported by the Australian Research Council.}

\date{Received 9 October 2000.}

\maketitle

Mackey's imprimitivity theorem characterises the unitary representations
of a locally compact group $G$ which have been induced from
representations of a closed subgroup $K$; Rieffel's influential
reformulation says that the group $C^*$-algebra $C^*(K)$ is Morita
equivalent to the crossed product $C_0(G/K)\times G$ \cite{rie74}. There
have since been many important generalisations of this theorem,
especially by Rieffel \cite{rie76, rie-applications} and by Green
\cite{green77, green-twisted}. These are all special cases of the
symmetric imprimitivity theorem of \cite{rae}, which gives a Morita
equivalence between two crossed products of induced $C^*$-algebras.

Quigg and Spielberg proved in \cite{qs}, by ingenious
but indirect methods,  that the symmetric imprimitivity
theorem, and hence all the other generalisations of Rieffel's
imprimitivity theorem, have analogues for reduced crossed products. A different Morita equivalence
between the same reduced crossed products was obtained by Kasparov \cite[Theorem~3.15]{kas}.

Here we
identify the representations which induce to regular representations under the Morita equivalence
of the symmetric imprimitivity theorem (see Theorem~\ref{ind=reg} and
Corollary~\ref{converse}), and thus obtain a direct proof of the theorem
of Quigg and Spielberg (see Corollary~\ref{firstcor}).  We discovered
Theorem~\ref{ind=reg} while trying to understand why Rieffel's theory of
proper actions in
\cite{rie-pr} gives an equivalence involving reduced crossed products
rather than full ones.

Theorem~\ref{ind=reg} has several other interesting applications. We can use
it to see, albeit somewhat indirectly, that regular representations
themselves nearly always induce to regular representations (see
Corollary~\ref{regtoreg}), and it gives a new proof of the main
theorem of \cite[\S4]{hrw} which avoids a complicated argument involving a
composition of crossed-product Morita equivalences (see
Corollary~\ref{corthird}). It also sheds light on constructions in
\cite{kw} and \cite{ER}, which, in various special cases of the symmetric
imprimitivity theorem, yield pairs of regular representations which induce to
each other (see Remarks~\ref{others}).

\subsection*{Notation} 
We shall use left Haar measure, because we always
do, and the right-regular representation $\rho$, because this is what seems
to come naturally out of our constructions; this means that
modular functions intrude, so that
$\rho=\rho^G:G\to U(L^2(G))$ is defined by
$(\rho_t\xi)(s)=\xi(st)\Delta_G(t)^{1/2}$. We also write $\rho$ for the representation on
$L^2(G,\H)$ defined by the same formula. With our conventions, the regular representation of
$A\times_\alpha G$ induced from a representation $\pi:A\to B(\H)$ is
the integrated form $\widetilde\pi\times \rho$ of the
covariant representation $(\widetilde\pi,\rho)$ of $(A,G,\alpha)$ on
$L^2(G,\H)$ in which
$(\widetilde\pi(a)\xi)(s)=\pi(\alpha_s(a))(\xi(s))$. 
Given a nondegenerate representation $\mu$ of $A$ we denote by $\bar\mu$ its extension 
to the multiplier algebra $M(A)$ of $A$.
We shall denote by
$\lt$ and $\rt$ any actions of groups by, respectively, left 
and right translation; thus if $H$ acts on the left of a locally
compact space $P$ and
$K$ acts on the right, we have $\lt_s(f)(p)=f(s^{-1}p)$ and
$\rt_t(f)(p)=f(pt)$ for
$f\in C_0(P)$.

\section{The Main Theorem}

We start with commuting free and proper actions of two locally compact
groups $H$ and
$K$ on the left and right of a locally compact space $P$, and commuting
actions
$\tau:H\to\Aut C$ and $\sigma:K\to\Aut C$ on a $C^*$-algebra $C$.
The \emph{induced $C^*$-algebra} $\Ind (C,\sigma)$ is the $C^*$-subalgebra
of $C_b(P,C)$ consisting of the functions $f$ such that
$f(pt) = \sigma_t^{-1}(f(p))$
for all $t\in K$ and $p\in P$, and such that the function 
$pK\mapsto\|f(p)\|$ vanishes at infinity on $P/K$. The
diagonal action $\lt\otimes\tau$ on $C_b(P,C)\subset M(C_0(P,C))$
restricts to a well-defined strongly continuous action of $H$ on $\Ind
(C,\sigma)$, which is characterised by
$(\lt\otimes\tau)_s(f)(p)=\tau_s(f(s^{-1}p))$. (The continuity of
this action was established in \cite[Lemma~A.1]{hrw}.)   Likewise,
$\Ind(C,\tau)$ consists of the bounded continuous functions $f:P\to
C$ such that
$f(sp) = \tau_s(f(p))$
for $s\in H$ and
$Hp\mapsto \|f(p)\|$ vanishes at infinity on $H\backslash P$, and we have
a natural action $\rt\otimes \sigma$ of $K$ on $\Ind (C,\tau)$ given by
$(\rt\otimes\sigma)_t(f)(p)=\sigma_t(f(pt))$.

The symmetric imprimitivity theorem of \cite[Theorem~1.1]{rae} shows
how to make $X_0:=C_c(P,C)$ into a pre-imprimitivity bimodule whose
completion $X$ implements a Morita equivalence between
$\Ind(C,\sigma)\times_{\lt\otimes\tau}H$ and
$\Ind(C,\tau)\times_{\rt\otimes\sigma}K$. Since a shortage of Greek
letters and left-right ambivalence have previously led to conflicts of
notation, it is worthwhile to record the formulas we
use:
\begin{align*}
a\cdot x(p) &= \int_H a(s,p) \tau_s(x(s^{-1}p))\,
\Delta_H(s)^{1/2}\, ds\\
x\cdot d(p) &=\int_K \sigma_t\big(x(pt) d(t^{-1},pt)\big)\,
\Delta_K(t)^{-1/2}\, dt\\
{}_{\Ind(C,\sigma)\times H}\langle{x},{y}\rangle(s,p) &= 
\int_K \sigma_t\big(x(pt)
\tau_s(y(s^{-1}pt)^*)\big)\, dt\,
\Delta_H(s)^{-1/2}\\
\langle{x},{y}\rangle_{\Ind(C,\tau)\times K}(t,p) &=
\int_H \tau_s\big(x(s^{-1}p)^*
\sigma_t(y(s^{-1}pt))\big)\,
ds\, \Delta_K(t)^{-1/2},
\end{align*}
where
$a\in C_c(H,\Ind(C,\sigma))\subset\Ind(C,\sigma)\times H$, $d\in
C_c(K,\Ind(C,\tau))\subset\Ind(C,\tau)\times K$, and $x,y$ belong to
$X_0=C_c(P,C)$.

We shall also need the one-sided version of this bimodule which is based
on the same space $Z_0:=C_c(P,C)$ but omits all mention of the group $H$; this
bimodule is the dual of the bimodule first considered in
\cite{rw85}. Thus we denote by $Z$ the completion of 
$Z_0=C_c(P,C)$ as a 
$\Ind(C,\sigma)$--$(C_0(P,C)\times_{\rt\otimes\sigma}K)$ imprimitivity 
bimodule.  We use exactly the same conventions as above, so that, for example,
\[
\langle{x},{y}\rangle_{C_0(P,C)\times K}(t,p) =
x(p)^*
\sigma_t(y(pt))\Delta_K(t)^{-1/2}.
\]

\begin{thm}\label{ind=reg}
Suppose $(\mu,U)$ is a covariant representation of
$\big(C_0(P,C),K,\rt\otimes\sigma\big)$ on $\H$. Then the representation
\[
X\dashind_{\Ind(C,\tau)\times
K}^{\Ind(C,\sigma)\times H}(\bar\mu|_{\Ind(C,\tau)}\times
U):\Ind(C,\sigma)\times_{\lt\otimes\tau}H\to
B(X\otimes_{\Ind(C,\tau)\times K}\H)
\]
is unitarily
equivalent to the right-regular representation 
\[
\big(Z\dashind_{C_0(P,C)\times
K}^{\Ind(C,\sigma)}(\mu\times U)\big)\widetilde{\;}\times\rho:
\Ind(C,\sigma)\times_{\lt\otimes\tau}H\to
B\big(L^2(H,Z\otimes_{C_0(P,C) \times K}\H)\big).
\]
\end{thm}

\begin{proof}
We shall write
$D:=\Ind(C,\tau)\times_{\rt\otimes\sigma}K$ and
$B:=C_0(P,C)\times_{\rt\otimes\sigma}K$; we are going to view
$X_0:=C_c(P,C)$ as the pre-Hilbert
$C_c(K,\Ind(C,\tau))$-module described above, and also as a dense subspace of the Hilbert
$B$-module
$Z$.

For $x\otimes_D h\in X_0\odot \H$, we define $W(x\otimes_D h):H\to
Z\otimes_B\H$ by
\[
W(x\otimes_D h)(s):=(\lt\otimes\tau)_s(x)\otimes_B h.
\]
We shall prove that $W$ extends to a unitary operator of
$X\otimes_{D}\H$ onto $L^2(H,Z\otimes_{B}\H)$ which intertwines the
given representations. 

We first prove that $W$ is well-defined and isometric: for both, it
suffices to show that
\begin{equation}\label{Wisometric}
\big(W(x\otimes_D h)\,\big|\,W(y\otimes_D
\bar\mu(f)k)\big)=\big(x\otimes_D h\,\big|\,y\otimes_D \bar\mu(f)k\big)
\end{equation}
for $x,y\in X_0$, $h,k\in \H$ and $f\in C_c(P)\subset M(C_0(P,C))$.
(Inserting the function $f$ allows us to deduce from the properness of the
actions that the integrands in the following calculations have compact
support.) To prove (\ref{Wisometric}), we note that
\begin{align*}
\big(x\otimes_B
h\,\big|\,y\otimes_B \bar\mu(f)k\big)&=\big(\mu\times U(\langle
y,x\rangle_B)h\,\big|\,\bar\mu(f)k\big)\\
&=\int_K\big(\mu(\langle
y,x\rangle_B(t))U_th\,|\,\bar\mu(f)k\big)\,dt\\
&=\int_K\big(\mu\big(p\mapsto
\overline{f(p)}y(p)^*\sigma_t(x(pt))
\Delta_K(t)^{-1/2}\big)U_th\,|\,k\big)\,dt.
\end{align*}
From this and an application of Fubini's
Theorem, we deduce that
\begin{align*}
\big(W(x\otimes_D h&)\,\big|\,W(y\otimes_D \bar\mu(f)k)\big)
=\int_H\big((\lt\otimes\tau)_s(x)\otimes_B
h\,\big|\,(\lt\otimes\tau)_s(y)\otimes_B \bar\mu(f)k\big)\,ds\\
&=\int_H\!\int_K\Big(\mu\big(p\mapsto
\overline{f(p)}\tau_s(y(s^{-1}p)^*
\sigma_t(x(s^{-1}pt)))\Delta_K(t)^{-1/2}\big)U_th\;
\Big|\;k\Big)\,dt\,ds\\
&=\int_K\big(\bar\mu(\overline{f}\langle
y,x\rangle_D(t))U_th\,\big|\,k\big)\,dt\\
&=\big(x\otimes_D
h\,\big|\,y\otimes_D \bar\mu(f)k\big),
\end{align*}
as desired.

We next prove that $W$ is surjective. We begin by observing that, with the
pointwise action $w\cdot b(s):=w(s)\cdot b$ and the inner product
\[
\langle v,w\rangle_B:=\int_H\langle v(s),w(s)\rangle_B\,ds,
\]
$C_c(H,Z_0)$ is a pre-Hilbert $B$-module; we denote its completion by
$L^2(H,Z)$. (Although we shall not need this, it might help to observe
that
$L^2(H,Z)$ is naturally isomorphic to the external Hilbert-module tensor
product
$L^2(H)_{\C}\otimes Z_B$ of, for example, \cite[\S3.4]{tfb}.) Writing out
the formulas shows that the map
$w\otimes_Bh\mapsto (s\mapsto w(s)\otimes_B h)$ is inner-product
preserving from
$L^2(H,Z)\otimes_B\H$ to $L^2(H,Z\otimes_B\H)$, and we can see by
considering elementary tensors
$w=\xi\otimes z\in C_c(H)\odot Z_0$ that it has dense range. So we can
view $L^2(H,Z\otimes_B\H)$ as $L^2(H,Z)\otimes_B\H$. Functions
$f\in C_c(P)\subset C_0(P)$ act as multipliers on $C_0(P,C)$ and hence
also on
$B=C_0(P,C)\times K$ and $X_0$, where the action is given by
$(x\cdot f)(t,p)=x(t,p)f(p)$. Since
\[
W(x\otimes_D \bar\mu(f)h)(s):=(\lt\otimes\tau)_s(x)\otimes_B \bar\mu(f)h
=((\lt\otimes\tau)_s(x)\cdot f)\otimes_B h,
\]
and since an isometry with dense range is surjective, it suffices for us
to show that
\[
L_0:=\sp\{(s,p)\mapsto \tau_s(x(s^{-1}p))f(p):x\in X_0,f\in C_c(P)\}
\]
is dense in $L^2(H,Z_B)$. 
If $v\in C_c(H\times P,C)\subset C_c(H,X_0)$
has $\supp v\subset M\times N$, then using the $L^1$-norm to estimate the
$C^*$-norm gives
\begin{align*}
\|v\|_B^2&=\int_H\|\langle v(s,\cdot),v(s,\cdot)\rangle_B\|\,ds\\
&\leq \int_H\Big(\int_K\sup_{p\in
P}\|v(s,p)^*\sigma_t(v(s,pt))\|\Delta_K(t)^{-1/2}\,dt\Big)\,ds\\
&\leq\mu_H(M)\int_{\{t\in K:N\cap
Nt^{-1}\not=\emptyset\}}\|v\|_\infty^2\Delta_K(t)^{-1/2}\,dt.
\end{align*}
Thus it suffices to show that we can approximate $v\in C_c(H\times
P,C)$ uniformly on a compact neighbourhood of $\supp
v$ by functions in
$L_0$ with support in that neighbourhood.

Because $H$ acts freely and properly on $P$, the map $(s,p)\mapsto
(sp,p)$ is a homeomorphism of $H\times P$ onto $P\times_{H\backslash
P}P:=\{(q,p)\in P\times P:Hq=Hp\}$. (The inverse is given by
$(q,p)\mapsto (\tr(q,p),p)$, where $\tr:P\times_{H\backslash
P}P\to H$ is the translation function
characterised by $q=\tr(q,p)p$; a
routine compactness argument shows that $\tr$ is continuous.) Thus the map
$\Phi:C_c(P\times_{H\backslash P}P,C)\to C_c(H\times P, C)$ defined by
$\Phi(w)(s,p)=\tau_s(w(s^{-1}p,p))$ is a  linear isomorphism which
preserves the kind of approximation we want. For $v\in C_c(H\times P,C)$,
we choose an extension $w$ of $\Phi^{-1}(v)$ to a function of a compact
support on $P\times P$, and now standard arguments show that we can
approximate
$w$ in
$C_c(P\times P, C)$ by functions $\sum x_i\otimes f_i$ in $C_c(P,C)\odot
C_c(P)$; then
$\Phi(\sum x_i\otimes f_i)$ is the required approximation to $v$. Thus
$W$ is surjective.  

It remains to check that $W$ intertwines the given representations as
claimed. Let
$a\in C_c(H,\Ind(C,\sigma))$, and for this calculation denote the
action $\lt\otimes\tau$ of $H$ on $\Ind(C,\sigma)$ by $\alpha$. Then
for $x\otimes_D h\in X_0\odot \H$, we have
\begin{equation}\label{WInd}
W\big(X\dashind (\bar\mu\times U)(a)(x\otimes_D
h)\big)(s)=W((a\cdot x)\otimes_D h)(s)=\alpha_s(a\cdot x)\otimes_B h.
\end{equation} 
On the other hand,
\begin{align}
\big(Z\dashind(\mu\times
U)\widetilde{\;}\times&\rho(a)\big)\big(W(x\otimes_D
h)\big)(s)\notag\\ 
&=\int_H Z\dashind(\mu\times
U)(\alpha_s(a(r)))\big(\rho_r(W(x\otimes_D h))(s)\big)\,dr\notag\\
&=\int_H Z\dashind(\mu\times
U)(\alpha_s(a(r)))\big(W(x\otimes_D
h)(sr)\Delta_H(r)^{1/2}\big)\,dr\notag\\
&=\int_H\big((\alpha_s(a(r))\cdot\alpha_{sr}(x))\otimes_B
h\big)\Delta_H(r)^{1/2}\,dr.\label{IndW}
\end{align}
Since we have
\begin{align*}
\alpha_s(a\cdot x)(p)&=\tau_s(a\cdot x(s^{-1}p))\\
&=\tau_s\Big(\int_H
a(r,s^{-1}p)\tau_r(x(r^{-1}s^{-1}p))\Delta_H(r)^{1/2}\,dr\Big)\\
&=\int_H
\big(\alpha_s(a(r))\cdot\alpha_{sr}(x)\big)(p)\Delta_H(r)^{1/2}\,dr,
\end{align*}
the only difference between (\ref{WInd}) and (\ref{IndW}) is the
location of the integral with respect to $\otimes_B h$. But the
integrands in both formulas are continuous and compactly supported, so
there is no difficulty verifying that they have the same inner product
with every vector of the form $y\otimes_B k\in Z_0\odot \H$, and hence
must be equal. We deduce that
\[
W\big(X\dashind (\bar\mu\times U)(a)\big)=\big(Z\dashind(\mu\times
U)\widetilde{\;}\times\rho\big)(a)W,
\]
and we have proved the Theorem.
\end{proof}

\begin{remark}
The referee pointed out that Theorem~\ref{ind=reg} suggests that there is, and would
follow from, an isomorphism
\begin{equation}\label{commdiag}
X\otimes_{\Ind(C,\tau)\times K}(C_0(P,C)\times K)
\cong Y\otimes_{\Ind(C,\sigma)} Z
\end{equation}
of right-Hilbert $(\Ind(C,\sigma)\times H)$--$(C_0(P,C)\times K)$ bimodules, where
$Y$ is the Hilbert bimodule of Green which induces representations of
$\Ind(C,\sigma)$ to regular representations of $\Ind(C,\sigma)\times H$. Such
isomorphisms have proved to be a powerful tool for studying the duality between
induction and restriction of representations \cite{kqr}. If applications arise which
require functorial properties of the equivalence in Theorem~\ref{ind=reg}, then
establishing such an isomorphism might be an efficient way to proceed; for our present
applications, Theorem~\ref{ind=reg} suffices.

\end{remark}

\section{Applications}

\subsection{The Theorem of Quigg and Spielberg} We retain the notation of
the preceding section. We shall say that a dynamical system $(A,G,\alpha)$ is
\emph{amenable} if the regular representation
induced from a faithful representation of $A$ is faithful on
$A\times_{\alpha}G$.

\begin{cor}[Quigg-Spielberg~\cite{qs}]\label{firstcor}
Denote by $I$ and $J$ the kernels of the quotient maps of 
$\Ind(C,\tau)\times_{\rt\otimes\sigma}K$ and
$\Ind(C,\sigma)\times_{\lt\otimes\tau}H$ onto the reduced crossed
products. Then
$X\dashind I=J$. In particular, this implies that
$\Ind(C,\tau)\times_{\rt\otimes\sigma,\; r}K$ is Morita equivalent to
$\Ind(C,\sigma)\times_{\lt\otimes\tau,\;r}H$, and that the system
$\big(\Ind (C,\sigma),H,\lt\otimes\tau\big)$ is amenable if and only if
$\big(\Ind (C,\tau),K,\rt\otimes\sigma\big)$ is amenable. 
\end{cor}

\begin{proof}
Every regular representation $\widetilde\pi\times\rho$ induced from
a faithful representation $\pi$ of $\Ind(C,\tau)$ has the same kernel
$I$, and
$X\dashind(\widetilde\pi\times\rho)$ has kernel $X\dashind I$. We
\emph{choose} $\pi$ to be the restriction $\bar\nu|_{\Ind(C,\tau)}$ of a
faithful nondegenerate representation of $C_0(P,C)$. Then
$(\widetilde\pi,\rho)$ is the restriction of the regular
representation $(\widetilde\nu,\rho)$ of
$\big(C_0(P,C),K,\rt\otimes\sigma\big)$, and we can apply
Theorem~\ref{ind=reg} with $(\mu,U)=(\widetilde\nu,\rho)$. We deduce
that, for this $\pi$, $X\dashind(\widetilde\pi\times\rho)$ is
equivalent to the right-regular representation of
$\Ind(C,\sigma)\times_{\lt\otimes\tau}H$ induced from the
representation $Z\dashind(\widetilde\nu\times\rho)$ of
$\Ind(C,\sigma)$. Thus
\begin{equation}\label{eq-ideals}
X\dashind I=\ker\big(X\dashind(\widetilde\pi\times\rho)\big)
=\ker\big(Z\dashind(\widetilde\nu\times\rho)\widetilde{\;}\times\rho\big)
\subset J.
\end{equation}
(If we knew that $Z\dashind(\widetilde\nu\times\rho)$ is faithful then we would have equality in \eqref{eq-ideals} above; instead we show that equality holds and deduce that $Z\dashind(\widetilde\nu\times\rho)$ is faithful.)
By symmetry
$\widetilde X\dashind J\subset I$, and now applying $X\dashind I$ we see that $J\subset X\dashind I$.

It follows from standard properties of the Rieffel correspondence
\cite[Proposition~3.25]{tfb} that the reduced crossed products are Morita
equivalent. Finally, if $\big(\Ind (C,\tau),K,\rt\otimes\sigma\big)$ is
amenable, then $I=\{0\}$, so $J=\{0\}$ and the system $\big(\Ind
(C,\sigma),H,\lt\otimes\tau\big)$ is amenable; the last part follows
by symmetry.
\end{proof}

One special case is worth mentioning
because the possibility of such a result was specifically mooted in
\cite[page~171]{rie-pr}, and because the proof we have given is more
direct than others.

\begin{cor}
Suppose $H$ acts freely and properly on $P$ and $\tau$ is an action of $H$
on a $C^*$-algebra $C$. Then
\[
C_0(P,C)\times_{\lt\otimes\tau}H=C_0(P,C)\times_{\lt\otimes\tau,\;r}H.
\]
\end{cor}

\begin{proof}
This is the special case of Corollary~\ref{firstcor} in which $K=\{e\}$;
the dynamical system $\big(\Ind (C,\tau),K,\rt\otimes\sigma\big)$ is then
trivially amenable.
\end{proof}

\subsection{Inducing regular representations} 

Because the symmetric
imprimitivity theorem passes to reduced crossed products, it is natural to ask
whether the symmetric imprimitivity theorem matches up the regular
representations themselves. More precisely, if
$\tilde\pi\times\rho$ is a regular representation of
$\Ind(C,\tau)\times_{\rt\otimes\sigma}K$, is the induced representation
$X\dashind(\tilde\pi\times\rho)$ of
$\Ind(C,\sigma)\times_{\lt\otimes\tau}H$ regular?
We can  settle this question, though in a rather roundabout
fashion (see Corollary~\ref{regtoreg} below).

We begin our analysis of this question by observing that Theorem~\ref{ind=reg}
can be used to characterise the representations of
$\Ind(C,\tau)\times_{\rt\otimes\sigma}K$ which induce to regular
representations of $\Ind(C,\sigma)\times_{\lt\otimes\tau}H$.

\begin{cor}\label{indreg}
Let $(\nu,V)$ be a covariant representation of
$(\Ind(C,\tau),K,\rt\otimes\sigma)$ on $\H$. Then the representation
$X\dashind(\nu\times V)$ is regular if and only if there is a covariant
representation
$(\mu,U)$ of
$(C_0(P,C),K,\rt\otimes\sigma)$ on $\H$ such that $(\nu,V)$ is equivalent to
$(\bar\mu|_{\Ind(C,\tau)}, U)$.
\end{cor}

\begin{proof}
Theorem~\ref{ind=reg} immediately gives the ``if'' direction. So suppose
$X\dashind(\nu\times V)$ is equivalent to the regular representation
$\tilde\pi\times\rho$ for some representation $\pi$ of $\Ind(C,\sigma)$. Let
$\mu\times U:=\widetilde Z\dashind_{\Ind(C,\sigma)}^{C_0(P,C)\times K}\pi$,
and note that $Z\dashind(\mu\times U)$ is equivalent to $\pi$. 
Theorem~\ref{ind=reg} implies that $X\dashind(\bar\mu|_{\Ind(C,\tau)}\times
U)$ is equivalent to $\tilde\pi\times\rho$, and applying $\widetilde
X\dashind$ shows that $(\nu,V)$ is
equivalent to $(\bar\mu|_{\Ind(C,\tau)}, U)$.
\end{proof}

\begin{cor}\label{converse}
Let $\pi$ be a nondegenerate representation of $\Ind(C,\tau)$. Then there
is a covariant representation $(\mu,U)$ of $(C_0(P,C),H,\lt\otimes\tau)$
such that
\[
X\dashind_{\Ind(C,\tau)\times
K}^{\Ind(C,\sigma)\times H}(\widetilde\pi\times\rho)=
\bar\mu|_{\Ind(C,\sigma)}\times U.
\]
\end{cor}

\begin{proof}
This follows from Corollary~\ref{indreg} by taking advantage of the 
symmetry of the situation. Suppose that $\pi$ is a representation of
$\Ind(C,\sigma)$ instead. Let $(\nu, V)=\widetilde
X\dashind(\tilde\pi\times\rho)$ and note that  $X\dashind(\nu\times V)$ is
regular because it is equivalent to $\tilde\pi\times\rho$. By
Corollary~\ref{indreg}
$(\nu, V)$ is the restriction of a covariant representation $(\mu, U)$ of
$(C_0(P,C),K,\rt\otimes\sigma)$.
\end{proof}

\begin{cor}\label{regtoreg}
Suppose, in addition to our standard assumptions on ${}_HP_K$, that $P$ is a locally trivial
principal $H$-bundle. Then for every nondegenerate representation $\pi$ of
$\Ind(C,\tau)$, the induced representation
$X\dashind(\widetilde\pi\times\rho)$ of
$\Ind(C,\sigma)\times_{\lt\otimes\tau}H$ is regular.
\end{cor}

By Green's imprimitivity theorem \cite[Theorem~6]{green-twisted}, it suffices
to construct a nondegenerate representation $\nu$ of $C_0(H)$ on the
Hilbert space
$X\otimes_D \H_\pi$ of $X\dashind(\widetilde\pi\times\rho)$ which
commutes with the action of $\Ind(C,\sigma)$ and, together with the unitary
part of $X\dashind(\widetilde\pi\times\rho)$, is covariant for the action
$\lt:H\to \Aut C_0(H)$: the imprimitivity theorem then implies that
$X\dashind(\widetilde\pi\times\rho)$ is induced from the subgroup $\{e\}$ of
$H$, and hence is regular. Let
$(\mu,U)$ be the covariant representation of $(C_0(P,C),H,\lt\otimes\tau)$
from Corollary~\ref{converse}.  We shall use $\mu$ and the copies of $H$
inside the principal bundle $P$ to construct the required representation
$\nu$ of $C_0(H)$. We make this  precise in the following Lemma:

\begin{lemma}\label{decomp}
Suppose $P$ is a locally trivial principal $H$-bundle and $(\mu,U)$ is a
covariant representation of $(C_0(P,C),H,\lt\otimes\tau)$ on $\H$. Then
there is a $(\mu,U)$-invariant subspace $\H_1$ of $\H$ and a representation
$\nu_1:C_0(H)\to B(\H_1)$ such that $(\nu_1,U|_{\H_1})$ is a covariant
representation of $(C_0(H),H,\lt)$ and each $\nu_1(f)$ commutes with each
$\mu(g)|_{\H_1}$. 
\end{lemma}

\begin{proof}
 We can  use a partition of unity on
$H\backslash P$ to write every function in the dense subalgebra $C_c(P,C)$ as
a sum of functions supported on $H$-saturated open subsets of $P$ which are
trivial as
$H$-bundles.  Since
$\mu$ is nondegenerate and in particular nonzero, $\mu$ must be nonzero
on one of these sets. More formally, there is an
$H$-saturated open set $N$ such that there is a bundle
isomorphism $\phi:N\to (H\backslash N)\times H$, and such that $\mu$ is not
identically zero on
$I_N:=\{g\in C_0(P,C):g(p)=0\mbox{ for }p\notin N\}$.  Because $N$ is
$H$-saturated, $I_N$ is invariant, and thus
$\H_1:=\clsp\{\mu(g)h:g\in I_N, h\in \H\}$ is a nonzero
$(\mu,U)$-invariant subspace of $\H$.

We now let $\phi_2:N\to (H\backslash N)\times H\to H$ denote the composition
of $\phi$ with the projection onto the second factor $H$, and define $\iota:C_0(H)\to ZM(I_N)$ by 
\[
(\iota(f)g)(p)=
\begin{cases}
f(\phi_2(p))g(p)&\text{for $p\in N$}\\
0&\text{for $p\notin N$,}
\end{cases}
\]
where $g\in I_N$.
Because $\mu|_{\H_1}$ is nondegenerate on $I_N$ it extends to a
representation $\bar\mu_1$ of $M(I_N)$. Since $\iota$ is nondegenerate as a homomorphism into $M(I_N)$ it follows that $\nu_1:=\bar\mu_1\circ\iota$
is a nondegenerate
representation of $C_0(H)$ on $\H_1$ whose range commutes
with every $\mu(g)|_{\H_1}$. Since $\phi$ is $H$-equivariant, $\iota$
is equivariant for the action $\lt$ of $H$ on $C_0(H)$ and the action
$\lt\otimes\tau$ of $H$ on $I_N\subset C_0(P,C)$; thus the
covariance of $(\mu,U)$ implies that $(\nu_1,U|_{\H_1})$ is a covariant
representation of $(C_0(H),H,\lt)$.
\end{proof}

\begin{proof}[Proof of Corollary~\ref{regtoreg}]
From Lemma~\ref{decomp} and a Zorn's Lemma argument, we obtain a
decomposition $\H=\bigoplus\H_i$ into $(\mu,U)$-invariant subspaces, each of
which admits a suitable nondegenerate representation $\nu_i:C_0(H)\to
B(\H_i)$. Now we just take $\nu:=\bigoplus \nu_i$, and apply Green's
imprimitivity theorem as described above.
\end{proof}

\begin{remark}
The local triviality hypothesis in Corollary~\ref{regtoreg} is a minor one,
and is automatic if $H$ is a Lie group, for example. Indeed, because the
action of $H$ is free and proper, $P$ is locally trivial if and only if the
orbit map
$q:P\to H\backslash P$ admits local continuous cross-sections
\cite[Proposition 4.65]{tfb}, and a theorem of Palais says that $q$ always
admits such sections when $H$ is a Lie group \cite[\S4.1]{palais}. 
\end{remark}

\begin{remarks}\label{others}
 One situation in which the induced representation is naturally regular
is that considered by  Kirchberg and Wassermann in \cite{kw}. Suppose $K$
is a closed subgroup of a locally compact group
$G$, $P=G$,
$H=G$ and 
$\tau_s=\id$ for all $s\in H$. Then $f\mapsto f(e)$ is an isomorphism of
$\Ind(C,\id)$ onto $C$, and the symmetric imprimitivity theorem says that
$\Ind(C,\sigma)\times_{\lt\otimes\id} G$ is  Morita equivalent to $C\times_\sigma
K$ \cite[\S1.5]{rae}.
If $\pi:C\to B(\H)$ is a
nondegenerate representation, then it is proved in \cite[Proposition
3.2]{kw} that $X\dashind(\widetilde \pi\times\lambda^K)$ is equivalent to a
regular representation\footnote{We warn that $\widetilde\pi$ has different
meaning in the presence of $\lambda$ and $\rho$.}. 

To see that $X\dashind(\widetilde
\pi\times\lambda^K)$ is regular using Theorem~\ref{ind=reg}, we define 
\begin{equation}\label{defmu}
(\mu(f)\xi)(t)=\pi(\sigma_t(f(t)))(\xi(t)).
\end{equation}
Then $(\mu,\rho^K)$ is a covariant representation of
$(C_0(G,C),K,\rt\otimes\sigma)$, and Theorem~\ref{ind=reg} implies that the
representation
$X\dashind(\bar\mu|_{\Ind(C,\id)}\times\rho^K)$ is regular. The extension
$\bar\mu$   is given on $C_b(G,C)$ by the same formula (\ref{defmu}), and
hence the isomorphism $\Ind(C,\id)\cong C$ carries $\bar\mu|_{\Ind(C,\id)}$
into $\widetilde\pi$. Thus $X\dashind(\widetilde \pi\times\rho^K)$ is regular,
and so is the equivalent representation $X\dashind(\widetilde
\pi\times\lambda^K)$. (As it stands, though, Theorem~\ref{ind=reg} does not
give the twisted version of this result given in \cite{kw}.)

In \cite[\S1]{ER} Echterhoff and Raeburn consider the special case 
where $H$ and $K$ are subgroups of the same locally compact group $G$ and
$P=G$. They construct a pair of regular representations of the induced
systems $(\Ind(C,\sigma),H)$ and $(\Ind(C,\tau),K)$, and show that these
induce to each other via the equivalence of the symmetric imprimitivity
theorem \cite[Theorem 1.4]{ER}. We can use our Theorem~\ref{ind=reg} to see
directly that the induced representations are regular: in the notation of
\cite{ER}, just define $\mu: C_0(G,D)\to\B\big(\Ind_H^G (L^2(K\times
H,\H_{\rho_1}))\big)$ by
$(\mu(g)\xi)(s)(k,h)=\rho_1(\beta_h(g(sh)))(\xi(s,k,h))$, and then $\mu$ is a
representation of $C_0(G, D)$ such that the restriction of $\bar\mu$ to the
induced algebra is $\Ind_H^G(\rho_1\otimes 1\otimes 1)$.
\end{remarks}

\subsection{Amenability of actions on $C_0(P)$-algebras}
We shall now show that the main
theorem of \cite[\S4]{hrw}, which is a generalisation of Corollary~\ref{firstcor} to
actions on
$C_0(P)$-algebras, can be deduced from Theorem~\ref{ind=reg}. We
recall the set-up of \cite{hrw}. We take ${}_HP_K$ as usual, but instead
of an arbitrary $C^*$-algebra $C$, we fix a
$C_0(P)$-algebra $A$; this means that there is a nondegenerate
injection $\iota_A$ of $C_0(P)$ into $ZM(A)$ (we write $f\cdot a$ for
$\iota_A(f)a$). We insist that the actions $\tau:H\to \Aut A$ and
$\sigma:K\to \Aut A$ commute and satisfy
\begin{equation}\label{consistency}
\tau_s(f\cdot a)=\lt_s(f)\cdot\tau_s(a)\ \mbox{ and }\ 
\sigma_t(f\cdot a)=\rt_t(f)\cdot\sigma_t(a).
\end{equation}
It is proved in \cite[\S3]{hrw} that both $\tau$ and $\sigma$ are proper
in the sense of \cite{rie-pr}, that there are strongly continuous actions
$\bar\tau:H\to\Aut A^\sigma$ and
$\bar\sigma:K\to \Aut A^\tau$ on the generalised fixed-point algebras of
\cite{rie-pr}, and that
$A^\sigma\times_{\bar\tau} H$ is Morita equivalent to
$A^\tau\times_{\bar\sigma}K$.

The next Corollary is Theorem~4.5
of
\cite{hrw}. There the one-sided case was proved first, by realising the
bimodule of
\cite{rie-pr} for the reduced crossed product as a quotient of the one from
\cite{rw85}; the general case was then deduced from several applications of
the one-sided case and a theorem of Combes \cite{com}. The proof of
Corollary~\ref{corthird}, on the other hand, uses only the general
results of \cite[\S1--2]{rae} and Theorem~\ref{ind=reg}; the
techniques, and in particular the indirect applications of
Theorem~\ref{ind=reg}, may be of independent interest. 

\begin{cor}\label{corthird}
In the above set-up, the system $(A^\sigma,H,\bar\tau)$ is amenable if and
only if
$(A^\tau,K,\bar\sigma)$ is amenable.
\end{cor}

As in \cite[\S3]{hrw}, we view $A^\tau$ as a quotient of $\Ind(A,\tau)$
--- indeed, for our purposes we could define $A^\tau$ this way, and avoid
all reference to proper actions. To see how this works, recall that the
nondegenerate action of $C_0(P)$ on $A$ induces a continuous surjection
$q_A:\Prim A\to P$, which is characterised by 
\[
q_A(I)=p\Longleftrightarrow\{f\cdot a\in C_0(P)\cdot A:f(p)=0\}\subset I;
\]
the hypothesis (\ref{consistency}) implies that $q_A$ is $H$- and
$K$-equivariant. It is proved in \cite[\S3]{hrw} (see
\cite[Proposition~3.6]{hrw} and the end of the proof of
\cite[Theorem~3.1]{hrw}) that the map $f\otimes a\mapsto f\cdot a$
extends to $C_b(P,A)\subset M(C_0(P)\otimes A)$, and induces an
equivariant isomorphism of
$\big(\Ind(A,\tau)/I(\tau),K,\rt\otimes\sigma\big)$ onto
$(A^\tau,K,\bar\sigma)$, where
\[
I(\tau):=\{b\in \Ind(A,\tau):b(q_A(I))\in I\mbox{ for all $I\in \Prim
A$}\}.
\]
This ideal is $(\rt\otimes\sigma)$-invariant, and hence by
\cite[Proposition~12]{green-twisted} the crossed product
$I(\tau)\times_{\rt\otimes\sigma}K$ embeds naturally as an ideal in
$\Ind(A,\tau)\times_{\rt\otimes\sigma}K$ with quotient 
\[
(\Ind(A,\tau)/I(\tau))\times_{\rt\otimes\sigma} K\cong
A^\tau\times_{\bar\sigma} K.
\]

To apply our theorem, we need a representation $\nu$ of $C_0(P,A)$ such
that $\bar\nu|_{\Ind(A,\tau)}$ has kernel $I(\tau)$. We choose a
nondegenerate representation $\nu$ of $C_0(P,A)$ with
\[
\ker\nu=I_\Delta:=\{b\in C_0(P,A):b(q_A(I))\in I\mbox{ for all $I\in
\Prim A$}\}.
\]
Then for $b\in\Ind(A,\tau)$, we have
\begin{align*}
\bar\nu(b)=0
&\Longleftrightarrow \nu(bc)=0\mbox{ for all $c\in C_0(P,A)$}\\
&\Longleftrightarrow b(q_A(I))c(q_A(I))\in I\mbox{ for all $I\in
\Prim A$, $c\in C_0(P,A)$}\\
&\Longleftrightarrow b(q_A(I))a\in I\mbox{ for all $I\in
\Prim A$, $a\in A$}\\
&\Longleftrightarrow b(q_A(I))\in I\mbox{ for all $I\in
\Prim A$,}
\end{align*}  
so $\ker\bar\nu|_{\Ind(A,\tau)}=I(\tau)$, as we wanted.

The representation $\big(\bar\nu|_{\Ind(A,\tau)}\big)\widetilde{\;}$ is
the restriction of (the extension to $M(C_0(P,A))$ of) the representation
$\widetilde\nu$ of $C_0(P,A)$; we aim to apply Theorem~\ref{ind=reg} to
the covariant representation $(\widetilde\nu,\rho)$ of
$\big(C_0(P,A),K,\rt\otimes\sigma\big)$. For this to be useful, we need
to know that $\ker (\widetilde\nu\times\rho)$ is the ideal in
$C_0(P,A)\times K$ corresponding to the diagonal ideal $I_\Delta$ in
$C_0(P,A)$.

\begin{lemma}\label{kerreg}
With the above notation, $\ker
(\widetilde\nu\times\rho)=I_\Delta\times_{\rt\otimes\sigma}K$. 
\end{lemma}

\begin{proof}
To avoid having to write out the opposite version of
Theorem~\ref{ind=reg}, we instead prove the equivalent assertion that the
regular representation $(\widetilde \nu,\rho)$ of the system 
$(C_0(P,A),H,\lt\otimes\tau)$ satisfies
$\ker(\widetilde\nu\times\rho)=I_\Delta\times_{\lt\otimes\tau}H$. To do this,
we apply Theorem~\ref{ind=reg} with
$K$ absent. Then
$X_0$ is the bimodule $Y_0:={}_{C_0(P,A)\times H}C_c(P,A)_{\Ind(A,\tau)}$
of
\cite[Theorem~2.2]{rw85},
$Z_0$ is the trivial bimodule ${}_{C_0(P,A)}C_c(P,A)_{C_0(P,A)}$, and
Theorem~\ref{ind=reg} says that
\[
Y\dashind(\bar\nu|_{\Ind(A,\tau)})\sim(Z\dashind
\nu)\widetilde{\;}\times\rho=\widetilde\nu\times\rho.
\]
Since $K=\{e\}$, we have $I(\sigma)=I_\Delta$,
and \cite[Corollary~2.1]{rae} implies that
$Y\dashind I(\tau)=I_\Delta\times H$. Thus our choice of $\nu$ implies that
$\ker\bar\nu|_{\Ind(A,\tau)}=I(\tau)$, and
\[
\ker(\widetilde\nu\times\rho)=
\ker\big(Y\dashind(\bar\nu|_{\Ind(A,\tau)})\big)
=Y\dashind I(\tau)=I_\Delta\times_{\lt\otimes\tau}H,
\]
as required.
\end{proof}

\begin{proof}[Proof of Corollary~\ref{corthird}]
Suppose that $(A^\tau,K,\bar\sigma)$ is amenable. Applying 
Theorem~\ref{ind=reg} to $(\widetilde\nu,\rho)$ shows that
\begin{equation}\label{dagger}
\ker\big(X\dashind((\bar\nu|_{\Ind(A,\tau)})\widetilde{\;}
\times\rho)\big)=
\ker\big((Z\dashind(\widetilde\nu\times\rho))\widetilde{\;}
\times\rho\big).
\end{equation}
Since $\ker\bar\nu|_{\Ind(A,\tau)}=I(\tau)$, $\bar\nu|_{\Ind(A,\tau)}$
factors through a faithful representation $\kappa_1$ of
$A^\tau\cong\Ind(A,\tau)/I(\tau)$; the amenability of
$(A^\tau,K,\bar\sigma)$ implies that $\widetilde{\kappa}_1\times\rho$
is faithful, or, equivalently, that $\ker\big((\bar\nu|_{\Ind(A,\tau)})\widetilde{\;}
\times\rho\big)=I(\tau)\times K$. Corollary~2.1 of \cite{rae} says
that the Rieffel correspondence $X\dashind$ carries $I(\tau)\times K$
to $I(\sigma)\times H$. Thus
\begin{equation}\label{dagger2}
I(\sigma)\times H=X\dashind(I(\tau)\times K)
=\ker\big(X\dashind((\bar\nu|_{\Ind(A,\tau)})\widetilde{\;}
\times\rho)\big).
\end{equation}

On the other hand, another application of \cite[Corollary~2.1]{rae}, with
$H$ missing this time, shows that the Rieffel correspondence $Z\dashind$
takes $I_\Delta\times K$ to $I(\sigma)$. Thus Lemma~\ref{kerreg} gives
\[
\ker\big(Z\dashind(\widetilde\nu\times\rho)\big)
=Z\dashind(I_\Delta\times K)=I(\sigma).
\]
So $Z\dashind(\widetilde\nu\times\rho)$ factors through a faithful
representation $\kappa_2$ of $A^\sigma\cong \Ind(A,\sigma)/I(\sigma)$,
and the representation
$\big(Z\dashind(\widetilde\nu\times\rho)\big)\widetilde{\;}
\times\rho$ appearing in (\ref{dagger}) factors through the regular
representation $\widetilde\kappa_2\times \rho$. Since we know from
(\ref{dagger}) and (\ref{dagger2}) that
$\big(Z\dashind(\widetilde\nu\times\rho)\big)\widetilde{\;}
\times\rho$ has kernel $I(\sigma)\times H$, $\widetilde\kappa_2\times
\rho$ must be faithful on $A^\sigma\times H\cong
(\Ind(A,\sigma)\times H)/(I(\sigma)\times H)$. Thus
$(A^\sigma,H,\bar\tau)$ is amenable.

The result follows by symmetry.
\end{proof}

\end{document}